\documentclass[10pt]{IEEEtran}
\usepackage{graphicx,epstopdf,multirow,multicol,subfigure,booktabs}

\usepackage{amssymb,amsthm}
\usepackage{algorithm}
\usepackage[noend]{algpseudocode}
\usepackage{color}
\usepackage{pgfplotstable}
\usepackage{pgfplots}
\pgfplotsset{compat=newest}
\usepackage{float}
\usepackage{setspace}
\usepackage[cmex10]{amsmath}
\usepackage{bm,mathtools}
\usepackage{cite}
\usepackage{times}
\usepackage{upgreek}
\usepackage{enumitem}

\algnewcommand\algorithmicinput{\textbf{Input:}}
\algnewcommand\Input{\item[\algorithmicinput]}

\algnewcommand\algorithmicoutput{\textbf{Output:}}
\algnewcommand\Output{\item[\algorithmicoutput]}

\usepackage{hyperref}
\usepackage{booktabs}
\allowdisplaybreaks
\graphicspath{figures/}

\theoremstyle{remark}

\usepackage{float}
\floatstyle{ruled}
\newfloat{model}{thp}{lop}
\floatname{model}{Model}

\newcommand{\omitit}[1]{}

\begin{document}
\title{Relaxations of AC Maximal Load Delivery \\ for Severe Contingency Analysis}

\author{
 \IEEEauthorblockN{Carleton Coffrin, Russell Bent, Byron Tasseff, Kaarthik Sundar, Scott Backhaus} \\
\IEEEauthorblockA{Los Alamos National Laboratory, 
Los Alamos, New Mexico, USA\\
\{cjc,rbent,btasseff,kaarthik,backhaus\}@lanl.gov}}

\markboth{Journal of \LaTeX\ Class Files,~Vol.~14, No.~8, April~2018}%
{Coffrin \MakeLowercase{\textit{et al.}}: AC Maximal Load Delivery for Severe Contingency Analysis}

\maketitle

\begin{abstract}
This work considers the task of finding an AC-feasible operating point of a severely damaged transmission network while ensuring that a maximal amount of active power loads can be delivered.  This AC Maximal Load Delivery (AC-MLD) task is a nonconvex nonlinear optimization problem that is incredibly challenging to solve on large-scale transmission system datasets.  This work demonstrates that convex relaxations of the AC-MLD problem provide a reliable and scalable method for finding high-quality bounds on the amount of active power that can be delivered in the AC-MLD problem.  To demonstrate their effectiveness, the solution methods proposed in this work are rigorously evaluated on 1000 $N$-$k$ scenarios on seven power networks ranging in size from 70 to 6000 buses.  The most effective relaxation of the AC-MLD problem converges in less than 20 seconds on commodity computing hardware for all 7000 of the scenarios considered.
\end{abstract}

\begin{IEEEkeywords}
Nonlinear Optimization, Convex Optimization, AC Maximal Load Delivery, N-k Contingency, Power System Restoration
\end{IEEEkeywords}

\vspace{-0.5cm}
\section*{Nomenclature}

\begin{IEEEdescription}[\IEEEusemathlabelsep\IEEEsetlabelwidth{$Y^s = g^s-$}]
  \item [{$N$}]  - The set of nodes 
  \item [{$E$, $E^R$}]  - The set of {\em from} and {\em to} branches 
  \item [{$G$}]  - The set of generators
  \item [{$L$}]  - The set of loads
  \item [{$H$}]  - The set of bus shunts
  \item [{$\bm i$}] - Imaginary number constant
  \item [{$S = p+ \bm iq$}] - AC power
  \item [{$V = v \angle \theta$}]  - AC voltage
  \item [{$Y = g + \bm ib$}]  - Branch admittance
  \item [{$Y^s_{ij}, Y^s_{ji}$}] - Line charging on {\em from} and {\em to} branches
  \item [{$T = t \angle \theta^t$}]  - Transformer properties
  \item [{$Y^s$}]  - Bus shunt admittance
  \item [{$W$}]  - Product of two AC voltages
  \item [{$s^u$}] - Branch apparent power thermal limit
  \item [{$\theta^{\Delta l}, \theta^{\Delta u}$}] - Voltage angle difference limits
  \item [{$S^d$}] - AC power demand
  \item [{$S^g$}] - AC power generation
  \item [{$c_0,c_1,c_2$}] - Generation cost coefficients 
  \item [{$\omega$}] - Load prioritization parameter 
   \item [{$\Re(\cdot), \Im(\cdot)$}] - Real and imaginary parts of a complex num.
   \item [{$(\cdot)^*$, $|\cdot|$}] - Conjugate and magnitude of a complex num.
  \item [{$x^l, x^u$}] - Lower and upper bounds of $x$, respectively
  \item [{$\bm x$}] - A constant value
\end{IEEEdescription}

\section{Introduction}

Restoring a transmission system after natural disasters, such as hurricanes, floods, or ice storms, is a challenging task with significant consequences to both human and economic welfare.  Indeed, a comprehensive approach to transmission system restoration requires considering many factors, such as component repair prioritization, cold load pickups, generation dispatch, and standing phase angle requirements \cite{9780780353978,COFFRIN2015144,7038388}.  This work focuses on a core subproblem of transmission system restoration, namely the  AC Maximal Load Delivery (AC-MLD) problem.  Informally, the AC-MLD problem is defined as follows. Given a severely damaged network, for example a network in which hundreds of components are out of service, the AC-MLD problem consists of determining the maximum amount of active power loads that can be served in the damaged network subject to operating requirements, such as bus voltage limits, branch thermal limits, and generator capability limits.

The AC-MLD problem is considered to be an incredibly challenging problem to solve on realistic datasets \cite{1265164,6345338,COFFRIN2015144}. The key challenge arises because of the significant number of damaged components, which renders the normal operation setpoint of little assistance in establishing a new AC-feasible operating point for the damaged network.  This task has been described as ``maddeningly difficult'' to do by hand \cite{1265164}.  Furthermore, this work highlights that AC transmission systems with realistic components, such as bus shunts and line charging, present additional challenges for finding AC-feasible solutions to the AC-MLD problem.

A natural approach to mitigate the computational challenges of the AC-MLD problem is to develop a linearized active-power-only approximation of the problem.  However, previous studies have demonstrated that this approach has significant limitations in severe contingency cases \cite{COFFRIN2015144,6345338}. To address these shortcomings, various extensions of active-power-only approximations, such as the LPAC approximation \cite{LPAC_ijoc,COFFRIN2015144}, have been considered. However, to the best of our knowledge, no work has attempted to solve the complete nonconvex nonlinear AC-MLD problem in its entirety without approximation.

In this work, we leverage recent developments in nonlinear programming and convex relaxations of the AC power flow equations \cite{Hijazi2017,7540869,1664986,Bai2008383,6507355,6507352,moment_hierarchy} to investigate the AC-MLD problem.  We propose a more detailed formulation of the AC-MLD problem and consider three relaxations of this formulation.  The experimental results demonstrate (i) the challenges of scaling the AC-MLD problem to realistic networks with thousands of buses, and (ii) the success of the proposed convex relaxations in providing strong upper bounds to the AC-MLD problem in just a few seconds.  A key advantage of the proposed solution methodology over previous works is that it is guaranteed to converge in polynomial time and provides an optimistic bound on the maximum amount of active power that can be served in the damaged network considering all of the typical operating requirements, such as bus voltage limits, branch thermal limits, and generator capability limits.

The rest of the paper is organized as follows. Section \ref{sec:context} discusses power system restoration broadly to provide context and motivation for the AC-MLD problem.  The scope of the AC power flow model is defined and a number of the feasibility challenges inherent in this model are discussed in Section \ref{sec:ac_model}. The AC-MLD formulation is proposed and various relaxations are derived in Section \ref{sec:ac_mld}. A detailed computational evaluation of all proposed formulations is conducted on seven realistic power networks in Section \ref{sec:eval}, and concluding remarks are provided in Section \ref{sec:conclusion}.

\section{Uses of the Maximal Load Delivery Problem}
\label{sec:context}

\begin{figure}[t]
    \begin{center}
    \includegraphics[width=8.75cm]{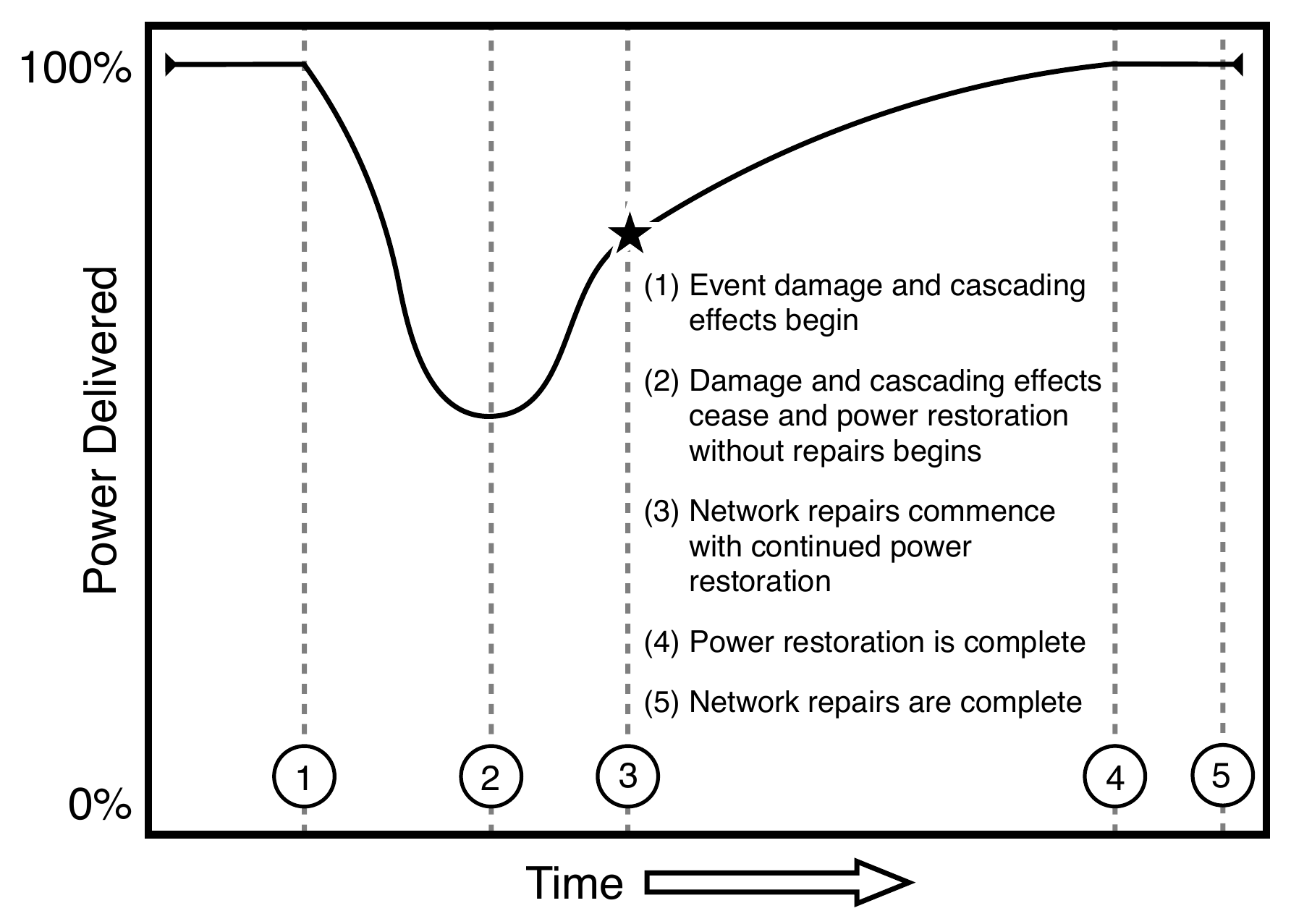}
    \end{center}
    \vspace{-0.3cm}
    \caption{A high-level illustration of a power system's response to an extreme event.  The star indicates the point in the damage and restoration process that this work seeks to quantify with an AC-MLD analysis.}
    \vspace{-0.3cm}
    \label{fig:restoration}
\end{figure}

The long-term goal of this line of work is to provide power system planning tools that can help network operators and policymakers understand and quantify how a power system will respond to extreme events, such as hurricanes, floods, and ice storms, where many components are out of service simultaneously.  A coarse overview of a typical power system response to an extreme event is presented in Figure \ref{fig:restoration}.  When the event begins, (1) power delivery decreases steadily as components are damaged and automatic projection devices lead to cascading effects.  At some point, these effects will cease, and the system will reach a stable operating point (2).  Some amount of power restoration can be executed without conducting repairs (3).  Then repairs are conducted in concert with power restoration until all of the desired power can be delivered (4).  Lastly, repairs will continue until all of the components in the system are operational (5).  Of course, the reality of an extreme event on a power system is more nuanced than the illustration in Figure \ref{fig:restoration} and features details such as discrete steps when breakers trip, nonlinear effects from cold load pickups, generation ramp rate limitations, and standing phase angle requirements for re-energizing transmission lines.

Addressing the complete scope of Figure \ref{fig:restoration} is a monumental task.  Consequently, research communities have emerged around each phase of the figure, including cascading failure analysis \cite{6112807,1428010,dobson2007complex}, power system restoration \cite{COFFRIN2015144,ThiebauxCHS13,7038388,5353669}, and restoration with repair crew routing \cite{van2011vehicle,coffrin2012last}.  This work's focus on the AC-MLD problem is most suitable to quantifying point (3) in Figure \ref{fig:restoration}.  Typically, a simulation tool will be used to generate a variety of component damage scenarios that may occur during an extreme event, and the AC-MLD problem is used to determine the maximal amount of load that can be delivered in these damage scenarios while accounting for the power restoration that can be executed without conducting repairs.  Additionally, the results of \cite{COFFRIN2015144,coffrin2012last,7038388} indicate that the AC-MLD problem is also a valuable building block for multi-timepoint power restoration problems that incorporate cold load pickups, generation ramp rates, and standing phase angle requirements.  The analysis of such extensions is left for future work.

\section{Scope of the AC Network Model}
\label{sec:ac_model}

In this work, we solve the AC-MLD problem at a level of network fidelity and scale that is comparable with commercial power system tools, such as PowerWorld, PSS/E, and PSLF.  However, in the interest of leveraging available network datasets, we adapt the AC Optimal Power Flow (AC-OPF) model of PowerModels \cite{pm_pscc}.  This formulation roughly corresponds to network models from PSS/E v33 and is an accurate mathematical model of the most common network components, namely generators, constant power loads, fixed bus shunts, $\pi$-model lines, and two-winding transformers.

\begin{model}[t]
\caption{The AC Optimal Power Flow (AC-OPF) Problem}
\label{model:ac_opf}
\begin{subequations}
\vspace{-0.2cm}
\begin{align}
\mbox{\bf variables: } & S^g_i (\forall i\in G), \; V_i (\forall i\in N)  \nonumber \\
%
\mbox{\bf minimize: } & \sum_{i \in G} \bm c_{2i} (\Re(S^g_i))^2 + \bm c_{1i}\Re(S^g_i) + \bm c_{0i} \label{ac_opf_obj} \\
\mbox{\bf subject to:} \nonumber
\end{align}
\vspace{-0.7cm}
\begin{align}
\phantom{12} & \bm {v^l}_i \leq |V_i| \leq \bm {v^u}_i \;\; \forall i \in N \label{ac_opf_1} \\
& \bm {S^{gl}}_i \leq S^g_i \leq \bm {S^{gu}}_i \;\; \forall i \in G \label{ac_opf_2}  \\
& \! \sum_{k \in G_{i}} \! S^g_k - \! \sum_{k \in L_{i}} \! \bm S^d_k - \! \sum_{k \in H_{i}} \bm Y^s_{k} |V_i|^2 = \!\!\!\!\!\!\!\!\!\! \sum_{\substack{(i,j)\in E_{i} \cup E^R_{i}}} \!\!\!\!\!\!\!\!\! S_{ij} \; \forall i\in N \label{ac_opf_3} \\ 
& S_{ij} = \left( \bm Y_{ij} + \bm Y^c_{ij} \right)^* \frac{|V_i|^2}{|\bm{T}_{ij}|^2} - \bm Y^*_{ij} \frac{V_i V_j^*}{\bm{T}_{ij}} \;\; (i,j)\in E \label{ac_opf_4}\\
& S_{ji} = \left( \bm Y_{ij} + \bm  Y^c_{ji} \right)^* |V_j|^2 - \bm Y^*_{ij} \frac{V_i^* V_j}{\bm{T}^*_{ij}} \;\; (i,j)\in E \label{ac_opf_5} \\
& |S_{ij}| \leq \bm {s^u}_{ij} \;\; \forall (i,j) \in E \cup E^R \label{ac_opf_6}  \\
& \bm {\theta^{\Delta l}}_{ij} \leq \angle (V_i V^*_j) \leq \bm {\theta^{\Delta u}}_{ij} \;\; \forall (i,j) \in E \label{ac_opf_7} 
\end{align}
\end{subequations}
\end{model}

The complete PowerModels AC-OPF formulation is presented in Model \ref{model:ac_opf}.  The objective function \eqref{ac_opf_obj} minimizes the active power generation costs.  Constraints \eqref{ac_opf_1} ensure that the bus voltage magnitudes are maintained within a desired operating range.  Constraints \eqref{ac_opf_2} provide generator operation limits for both active and reactive injection.   Constraints \eqref{ac_opf_3} capture power balance from Kirchhoff's current law, and constraints \eqref{ac_opf_4}-\eqref{ac_opf_5} model Ohm's law over branches.  Finally, constraints \eqref{ac_opf_6}-\eqref{ac_opf_7} ensure that branch thermal limits and voltage angle difference limits are enforced.  A detailed derivation of this formulation and its notations is available in \cite{7271127}.

It is important to emphasize that the AC modeling details, such as bus shunts (i.e., $\bm Y^s$ in \eqref{ac_opf_3}) and line charging (i.e., $\bm Y^c$ in \eqref{ac_opf_4}-\eqref{ac_opf_5}), present unique challenges for developing an AC-MLD formulation.  This is in significant contrast to AC-OPF formulations, where incorporating these terms is a trivial extension.

\subsection{AC Feasibility Challenges under Contingencies}
\label{sec:ac_model_probs}
At first glance, adapting Model \ref{model:ac_opf} to an AC-MLD formulation seems straightforward.  One can make the constant power loads flexible by introducing a decision variable $z^d_k \in (0,1)$ and updating the power balance constraint as follows:
\begin{align}
\!\!\!\!\! \sum_{k \in G_{i}} S^g_k - \sum_{k \in L_{i}} z^d_k {\bm S^d_k} - \sum_{k \in H_{i}} \bm Y^s_{k} |V_i|^2 = \!\!\!\!\!\!\!\!\! \sum_{\substack{(i,j)\in E_{i} \cup E^R_{i}}} \!\!\!\!\!\!\!\! S_{ij} \;\; \forall i\in N
\end{align}
This modification introduces bus-by-bus constant power factor loading into the model, and intuitively, one expects there to be a reduced loading that is sufficient to find a new AC operating point.  However, after considering this simple extension on several networks under thousands of contingencies, it was observed that this model is insufficient to ensure feasibility in all cases.  

\begin{figure}[t]
    \begin{center}
    \includegraphics[width=5.0cm]{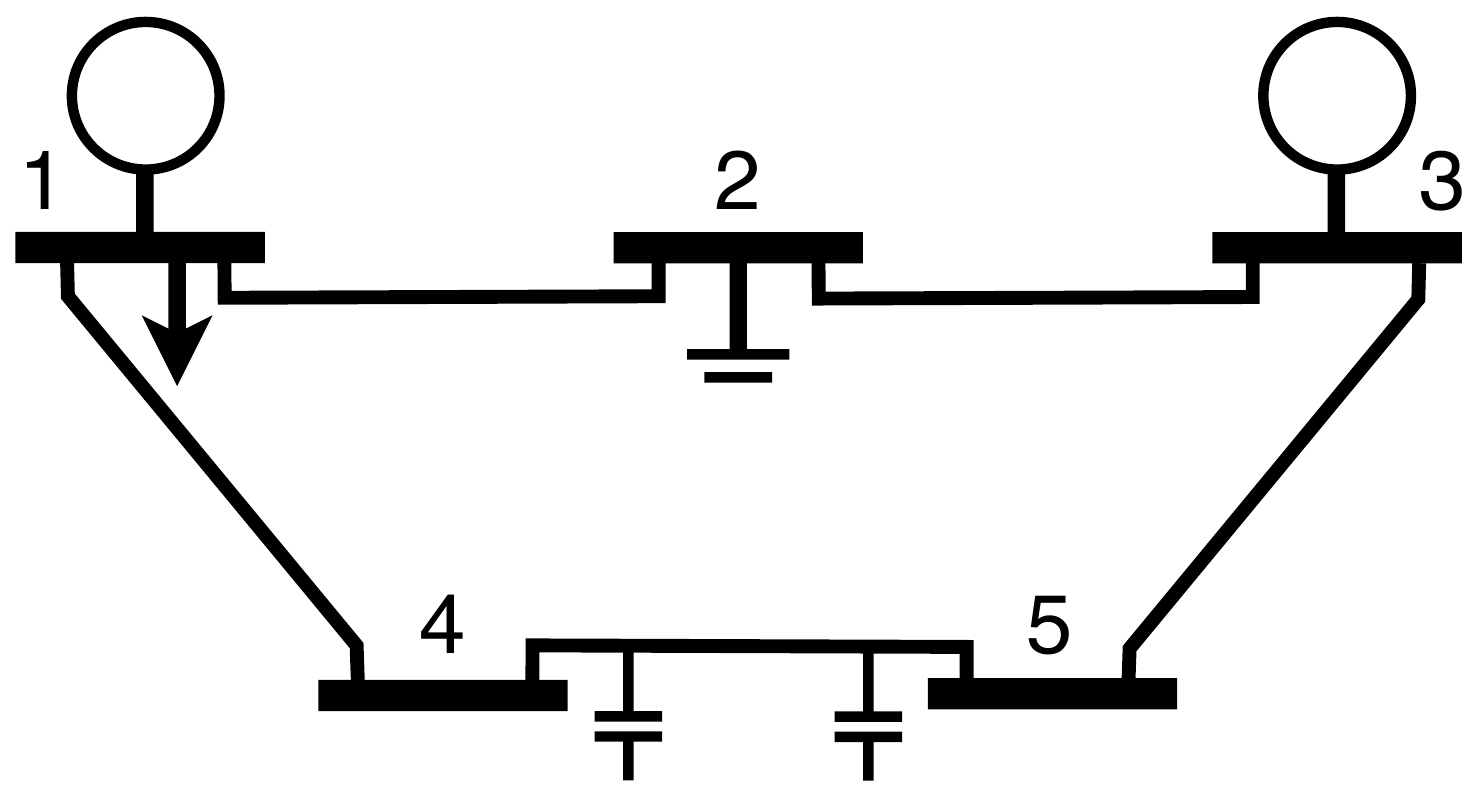}
    \end{center}
    \vspace{-0.3cm}
    \caption{A network diagram of a five-bus example illustrating AC feasibility challenges under contingencies.}
    \vspace{-0.2cm}
    \label{fig:5_bus}
\end{figure}

Figure \ref{fig:5_bus} presents a simple five-bus AC power network featuring bus shunts and line charging to illustrate three unexpected contingency cases that have been identified.  To make this example mathematically precise, we assume the following network data, which are representative of realistic network datasets:
\begin{subequations}
\begin{align}
& 0.9 \leq |V_i| \leq 1.1 \;\; \forall i \in \{1,2,3,4,5\}  \\
& 10 - \bm i \infty \leq S^g_3 \leq 100 + \bm i \infty \\
& \bm Y_{4,5} = 0 - \bm i 25, \bm Y^c_{4,5} = 0 + \bm i 0.04, \bm Y^c_{5,4} = 0 + \bm i 0.04
\end{align}
\end{subequations}

\paragraph*{Generator Injection Bounds}
On the five-bus example network, consider the outage contingency of line $2$-$3$ and line $3$-$5$. This isolates bus 3 and its associated generating unit. The problem in this contingency is that the generator has a non-zero active power injection bound.  Hence, when there is no load to serve, the constraints cannot be satisfied. Specifically, in this contingency, constraints \eqref{ac_opf_2} and \eqref{ac_opf_3} result in
\begin{subequations}
\begin{align}
& 10 \leq \Re(S^g_3) \leq 100 \\
& \Re(S^g_3) = 0
\end{align}
\end{subequations}
which has no feasible assignment.  This indicates that model feasibility requires all generator operating ranges to include the zero value or generators should be uncommitted from the dispatch scenario.

\paragraph*{Bus Shunt Shedding}
On the five-bus example network, consider the outage contingency of line $1$-$2$ and line $2$-$3$.  This isolates bus 2 and its associated fixed bus shunt. The key problem in this contingency is that a fixed shunt results in a nonzero active and reactive power injection at a bus. Hence when no generation is available, power balance cannot be satisfied. Specifically, in this contingency, constraints \eqref{ac_opf_1} and \eqref{ac_opf_3} result in
\begin{subequations}
\begin{align}
& 0.9 \leq |V_2| \leq 1.1 \\
& \bm Y^s_2 |V_2|^2 = 0
\end{align}
\end{subequations}
which has no feasible assignment.  This indicates that model feasibility requires the option of shedding fixed bus shunts from the network.

\paragraph*{Bus Voltage Bounds}
On the five-bus example network, consider the outage contingency of line $1$-$4$ and line $3$-$5$.  This isolates line $4$-$5$, which has a nonzero line charging value.  The key problem in this case is that line charging makes the Ohm's law constraints impossible to satisfy when no power flows over the line.  Specifically, in this contingency, constraints \eqref{ac_opf_1}, \eqref{ac_opf_4}, and \eqref{ac_opf_5} result in
\begin{subequations}
\begin{align}
& 0.9 \leq |V_4| \leq 1.1 \label{ac_lc_1} \\
& 0.9 \leq |V_5| \leq 1.1 \label{ac_lc_2} \\
& 0 = \left( \bm Y_{4,5} + \bm {Y^c}_{4,5} \right)^* |V_4|^2 - \bm Y^*_{4,5} V_4 V_5^* \label{ac_lc_3} \\
& 0 = \left( \bm Y_{4,5} + \bm {Y^c}_{5,4} \right)^* |V_5|^2 - \bm Y^*_{4,5} V_4^* V_5 \label{ac_lc_4}
\end{align}
%
An analysis of the active power aspect of \eqref{ac_lc_3}-\eqref{ac_lc_4} reveals that one of the following must hold:
%
\begin{align}
|V_4| = 0 \vee |V_5| = 0 \vee \angle(V_i V^*_j) = 0 \label{ac_lc_5}
\end{align}
%
while an analysis of the reactive power aspect of \eqref{ac_lc_3}-\eqref{ac_lc_4} reveals that both of the following must hold:
%
\begin{align}
|V_4| = |V_5| \wedge \angle(V_i V^*_j) \neq 0 \label{ac_lc_6}
\end{align}
\end{subequations}
Combining the requirements of \eqref{ac_lc_1}, \eqref{ac_lc_2}, \eqref{ac_lc_5}, and \eqref{ac_lc_6} demonstrates that this system of equations has no feasible assignment.  This indicates that model feasibility requires the option of either removing lines with charging values or relaxing the lower bounds of bus voltage magnitudes.  These observations indicate that the AC-MLD problem requires a disjunctive optimization formulation that is considerably more difficult to solve than traditional AC-OPF formulations.

\section{AC Maximal Load Delivery Formulations}
\label{sec:ac_mld}

Motivated by the observations in Section \ref{sec:ac_model}, this work introduces the following decision variables into the AC-MLD problem: 
\begin{subequations}
\begin{align}
z^v_i \in \{0, 1\} \;\; \forall i \in N \label{ac_mls_v1} \\
z^g_i \in \{0, 1\} \;\; \forall i \in G \label{ac_mls_v2} \\
z^d_i \in (0, 1) \;\; \forall i \in L \label{ac_mls_v3} \\
z^s_i \in (0, 1) \;\; \forall i \in H \label{ac_mls_v4}
\end{align}
\end{subequations}
In other words, discrete on/off variables for generators and bus voltages \eqref{ac_mls_v1}--\eqref{ac_mls_v2} and continuous shedding variables for loads and bus shunts \eqref{ac_mls_v3}--\eqref{ac_mls_v4}.
These additional variables allow the AC-MLD formulation to remove generators, bus shunts, and entire buses if needed to ensure model feasibility.  This additional flexibility resolves the problematic cases highlighted in Section \ref{sec:ac_model_probs}.

Additionally, the parameters $\bm \omega_i \geq 0 \; \forall i \in L$ are used to specify load restoration priorities.  If no load priorities are available, a setting of $\bm \omega_i = 1$ can be used to give all loads an equal priority.

\subsection{AC Maximal Load Delivery}

\begin{model}[t]
\caption{The AC Maximal Load Delivery (AC-MLD) Problem}
\label{model:ac_mls}
\begin{subequations}
\vspace{-0.2cm}
\begin{align}
\mbox{\bf variables: } & S^g_i (\forall i\in N), \; V_i (\forall i\in N)  \nonumber \\
& z^v_i \in \{0, 1\} (\forall i \in N),\; z^g_i \in \{0, 1\} (\forall i \in G) \nonumber \\
& z^d_i \in (0, 1) (\forall i \in L),\; z^s_i \in (0, 1) (\forall i \in H) \nonumber \\
%
\mbox{\bf maximize:} & \left( \sum_{i \in N} z^v_i, \sum_{i \in G} z^g_i, \sum_{i \in H} z^s_i, \sum_{i \in L} \bm \omega_i |\Re(\bm S^d_i)| z^d_i \right)  \label{ac_mls_obj} \\
\mbox{\bf subject to:} \nonumber
\end{align}
\vspace{-0.5cm}
\begin{align}
& z^v_i \bm {v^l}_i \leq |V_i| \leq z^v_i \bm {v^u}_i \;\; \forall i \in N \label{ac_mls_1} \\
& z^g_i \bm {S^{gl}}_i \leq S^g_i \leq z^g_i \bm {S^{gu}}_i \;\; \forall i \in G \label{ac_mls_2} \\
& \sum_{k \in G_{i}} S^g_k - \sum_{k \in L_{i}} z^d_k {\bm S^d_k} - \sum_{k \in H_{i}} z^s_k \bm Y^s_{k} |V_i|^2 = \!\!\!\!\!\!\!\!\! \sum_{\substack{(i,j)\in E_{i} \cup E^R_{i}}} \!\!\!\!\!\!\!\! S_{ij} \; \forall i \!\in\! N \label{ac_mls_3} \\
& \mbox{\eqref{ac_opf_4}--\eqref{ac_opf_7}} \nonumber
\end{align}
\end{subequations}
\end{model}

The complete AC-MLD formulation is presented in Model \ref{model:ac_mls}.  The multi-objective \cite{GoalProgramming1978} nature of \eqref{ac_mls_obj} strives to (i) keep as many buses in the network as possible, (ii) keep as many generators in the network as possible, (iii) remove the least amount of bus shunt admittance as possible, and (iv) remove the least amount of active power load as possible.  Constraints \eqref{ac_mls_1} ensure that bus voltage magnitudes are maintained within their operating range unless the bus is out of service, in which case the voltage magnitude is set to zero.  Constraints \eqref{ac_mls_2} ensure that generators are operated within their injection range unless the generator is out of service, in which case the generator outputs are set to zero.  Finally, constraints \eqref{ac_mls_3} incorporate continuous shunt and load removal into the power balance constraint.  The remaining constraints \eqref{ac_opf_4}--\eqref{ac_opf_7} are identical to the AC-OPF formulation and implement Ohm's law and the branch flow limits.

Model \ref{model:ac_mls} is a multi-objective nonconvex mixed-integer nonlinear program (MINLP), a class of problems that are typically NP-hard.  Hence, the existence of a scalable and reliable algorithm for solving Model \ref{model:ac_mls} is not obvious \cite{7063278}.  Developing a computationally efficient algorithm to solve Model \ref{model:ac_mls} involves addressing the following three aspects of the model: (i) the multi-objective function \eqref{ac_mls_obj}, (ii) the discrete variables (i.e., $z^v, z^g$), and (iii) the nonconvex constraints \eqref{ac_mls_3} and \eqref{ac_opf_4}--\eqref{ac_opf_7}.  The forthcoming sections discuss established approaches for addressing these model features. Each of these approaches increases the scalablity and reliability of solving the problem at the expense of omitting some aspect of the original AC-MLD formulation.

\subsection{Single-Objective AC Maximal Load Delivery} \label{subsec:single}
Although it is preferable to solve AC-MLD's multi-objective form and despite the recent algorithmic advances in multi-objective optimization \cite{ijoc.2015.0657,Boland2016}, there is only limited support for multi-objective functions in industrial optimization software.  To work around this practical limitation, we propose modeling the multi-objective function \eqref{ac_mls_obj} as a single weighted objective function as follows:
\begin{align}
\! \sum_{i \in N} \bm M^v z^v_i + \! \sum_{i \in G} \bm M^g z^g_i + \! \sum_{i \in H} \bm M^s z^s_i + \! \sum_{i \in L} \bm \omega_i |\Re(\bm S^d_i)| z^d_i \label{eq:single_obj}
\end{align}
The key challenge in using this representation is to find appropriate values of the scaling parameters $\bm M^v$, $\bm M^g$, and $\bm M^s$.  For the datasets considered in this work, we found the following method of setting these parameters to be suitable:
$\bm M^s = 10 \cdot \max_{i \in L}(\bm \omega_i |\Re(\bm S^d_i)|)$, $\bm M^g = \bm M^s$, and $\bm M^v = 10 \cdot \bm M^s$.
%
%
Intuitively, these weights strive to keep the bus, generator, and shunt settings at their nominal pre-contingency state and treat load removal as the primary control parameter.  It is useful to note that if numerical issues are encountered using this weighted objective function \eqref{eq:single_obj}, then the proposed lexicographic multi-objective \eqref{ac_mls_obj} can be solved by a sequence of single-objective optimization problems \cite{GoalProgramming1978}, which will increase numerical stability and forgo the weighting parameters.  Throughout the rest of this work, the AC-MLD formulation and its variants are assumed to be single-objective optimization problems.

\subsection{Relaxation of Discrete Variables}

To address the discrete variables (i.e., $z^v, z^g \in \{0, 1\}$) in the AC-MLD formulation in Model \ref{model:ac_mls}, the standard integrality relaxation is leveraged \cite{big_m}, namely
\begin{align}
& z^v_i, z^g_i \in \{0, 1\} \Rightarrow z^v_i, z^g_i \in (0, 1)
\end{align}
This relaxation allows the variables to take a continuous range of values between 0 and 1.  Throughout this work, {\em C} is used to annotate models where the discrete variables have been relaxed to continuous variables (e.g., AC-MLD-C).

\subsection{Relaxations of Nonconvex Constraints}

The sources of nonconvexity in Model \ref{model:ac_mls} take the form of bilinear products of continuous variables (e.g., $z^s|V|^2$, $V_i^* V_j$).  Recent developments in convex relaxation of the AC-OPF problem are utilized to address these nonconvexities \cite{7271127,7540869,1664986,Bai2008383,moment_hierarchy}.

\paragraph*{The Second-Order Cone Relaxation} 
In the interest of performance and scalablity, we choose to develop a model based on the Second-Order Cone (SOC) relaxation of the AC power flow equations \cite{1664986}.  Although some relaxations are stronger than SOC \cite{7271127,6345272,moment_hierarchy} and others are faster than SOC \cite{7540869}, the SOC relaxation is selected because it provides an appealing tradeoff between bounding strength and runtime performance.

The first insight of the SOC relaxation is that the voltage product terms $V_i^* V_j$ can be lifted into a higher dimensional $W$-space as follows:
\begin{subequations}
\begin{align}
|V_i|^2 &\Rightarrow W_{ii} \;\; \forall i \in N \label{w_lift_1} \\
V_i V^*_j &\Rightarrow W_{ij} \;\; \forall (i,j) \in E \label{w_lift_2}
\end{align}
\end{subequations}
Note that lifting Model \ref{model:ac_opf} into the $W$-space makes all of the nonconvex constraints linear.  
The second insight of the SOC relaxation is that this $W$-space relaxation can be strengthened by adding the valid inequality
\begin{subequations}
\begin{align}
& |W_{ij}|^2 \leq W_{ii}W_{jj} \label{w_soc}
\end{align}
\end{subequations}
which is a convex rotated SOC constraint that is supported by a wide variety of industrial-grade optimization tools.

After applying \eqref{w_lift_1}, \eqref{w_lift_2}, and \eqref{w_soc} to Model \ref{model:ac_mls}, one nonconvex term still remains in the power balance constraint because of the product of the shunt shedding variable $z^s_i$ and the lifted voltage magnitude variable $W_{ii}$ as follows:
\begin{align}
& \!\! \sum_{k \in G_{i}} S^g_k - \! \sum_{k \in L_{i}} z^d_k {\bm S^d_k} - \! \sum_{k \in H_{i}} z^s_k \bm Y^s_{k} W_{ii} = \!\!\!\!\!\!\!\!\! \sum_{\substack{(i,j)\in E_{i} \cup E^R_{i}}} \!\!\!\!\!\!\!\! S_{ij} \; \forall i \!\in\! N
\end{align}
Observing that the $z^s_i$ and $W_{ii}$ variables have tight bounds, McCormick inequalities \cite{MacC76} can be leveraged to provide a convex relaxation of this nonconvex variable product.  The general form of the McCormick envelope is given by
\begin{equation*}
\tag{M-CONV}
\langle xy \rangle^M \equiv
\begin{cases*}
xy  \geq  \bm {x^l}y + \bm {y^l}x - \bm {x^l}\bm {y^l}\\
xy  \geq  \bm {x^u}y + \bm {y^u}x - \bm {x^u}\bm {y^u}\\
xy  \leq  \bm {x^l}y + \bm {y^u}x - \bm {x^l}\bm {y^u}\\
xy  \leq  \bm {x^u}y + \bm {y^l}x - \bm {x^u}\bm {y^l}
\end{cases*}
\end{equation*}
In this case, a new variable $W^s_{k}$ is introduced to represent the product of $z^s_k$ and $W_{ii}$ as follows:
\begin{align}
& W^s_{k} =  \langle z^s_k W_{ii} \rangle^M \;\; \forall k \in L
\end{align}
With this final relaxation, a convex SOC relaxation of the AC-MLD problem is developed.

The complete SOC-MLD formulation is presented in Model \ref{model:soc_mls}.  Constraints \eqref{soc_mls_1} enforce the bus voltage magnitudes' operating range in the lifted $W$-space variables.  Constraints \eqref{soc_mls_2} provide the convex relaxation of the shunt shedding and voltage magnitude variable products, and constraints \eqref{soc_mls_3} incorporate the relaxation variable $W^s_{k}$ into the power balance constraints.  Constraints \eqref{soc_mls_4}--\eqref{soc_mls_7} implement the established $W$-space convex relaxation of \eqref{ac_opf_4}, \eqref{ac_opf_5}, and \eqref{ac_opf_7}. A detailed derivation of these constraints is available in \cite{7271127}.

\begin{model}[t]
\caption{The Maximal Load Delivery Relaxation (SOC-MLD)}
\label{model:soc_mls}
\begin{subequations}
\vspace{-0.2cm}
\begin{align}
\mbox{\bf variables: } & S^g_i (\forall i\in N), \; W_{ij} (\forall (i,j)\in E), \; W_{ii} (\forall i \in N) \nonumber \\
& W^s_{i} (\forall i \in L) \nonumber \\
& z^v_i \in \{0, 1\} (\forall i \in N),\; z^g_i \in \{0, 1\} (\forall i \in G) \nonumber \\
& z^d_i \in (0, 1) (\forall i \in L),\; z^s_i \in (0, 1) (\forall i \in H) \nonumber \\
%
\mbox{\bf maximize: } & \mbox{\eqref{ac_mls_obj}}
\nonumber \\
\mbox{\bf subject to:} \nonumber
\end{align}
\vspace{-0.7cm}
\begin{align}
& \mbox{\eqref{ac_opf_6}, \eqref{ac_mls_2}} \nonumber \\
& z^v_i (\bm {v^l}_i)^2 \leq W_{ii} \leq z^v_i (\bm {v^u}_i)^2 \;\; \forall i \in N \label{soc_mls_1} \\
& W^s_{k} = \langle z^s_k W_{ii} \rangle^M \;\; \forall k \in L_i,  \;  \forall i \in N \label{soc_mls_2} \\ 
& \!\! \sum_{k \in G_{i}} S^g_k - \!\! \sum_{k \in L_{i}} z^d_k {\bm S^d_k} - \!\! \sum_{k \in H_{i}} \bm Y^s_{k} W^s_{k} = \!\!\!\!\!\!\!\!\! \sum_{\substack{(i,j)\in E_{i} \cup E^R_{i}}} \!\!\!\!\!\!\!\! S_{ij} \;\; \forall i \!\in\! N \label{soc_mls_3} \\
& S_{ij} = \left( \bm Y_{ij} + \bm Y^c_{ij} \right)^* \frac{W_{ii}}{|\bm{T}_{ij}|^2} - \bm Y^*_{ij} \frac{W_{ij}}{\bm{T}_{ij}} \;\; (i,j)\in E \label{soc_mls_4}\\
& S_{ji} = \left( \bm Y_{ij} + \bm Y^c_{ji} \right)^* W_{jj} - \bm Y^*_{ij} \frac{W^*_{ij}}{\bm{T}^*_{ij}} \;\; (i,j)\in E \label{soc_mls_5} \\
& \tan(-\bm {\theta^{\Delta l}}_{ij}) \Re(W_{ij}) \leq \Im(W_{ij}) \leq \tan(\bm {\theta^{\Delta u}}_{ij}) \Re(W_{ij}) \label{soc_mls_6} \\
& |W_{ij}|^2 \leq W_{ii}W_{jj} \;\; \forall (i,j)\in E \label{soc_mls_7}
\end{align}
\end{subequations}
\end{model}

\subsection{Formulations Summary}

In summary, this work considers the following four variants of the proposed AC-MLD formulation from Model \ref{model:ac_mls}:
\begin{enumerate}[label=(\roman*)]
    \item {\em AC-MLD} - a nonconvex MINLP that is a challenging model to solve at scale.
    \item {\em AC-MLD-C} - a nonconvex nonlinear program (NLP) that is more likely to converge at scale but may not satisfy the integrality requirements of AC-MLD.
    \item {\em SOC-MLD} - a mixed-integer SOC program (MISOCP) that is guaranteed to converge but provides only an upper bound on the optimal value of AC-MLD.
    \item {\em SOC-MLD-C} - a SOC program (SOCP) that is a weaker but much more scalable variant of SOC-MLD.
\end{enumerate}
A core motivation for considering all of these variants is that they provide a range of tradeoffs in model accuracy and performance.  This allows us to clearly quantify how each relaxation affects the AC-MLD solution quality. 

\subsection{Data Processing}

Traditional AC-OPF studies are conducted on test cases that have been curated by subject matter experts to ensure that the network data is meaningful and of high quality \cite{nesta}.  However, in the context of severe contingency analysis, the AC-MLD is subject to network data where hundreds of the components are out of service, which can lead to unexpected network structures.  A notable example of this is the occurrence of {\em dangling buses}.  These buses are connected to the network by a single branch and have no load or generation and hence have limited utility.

In this work, we observed that conducting the following data processing steps on damaged networks greatly aided the convergence of all formulations considered herein. First, all dangling buses should be put out of service.  Second, connected components without load or generation should be put out of service.  Third, each connected component in the network should be solved independently. For simplicity, in this work we solved only the largest connected component.  

\section{Computational Evaluation}
\label{sec:eval}

To evaluate the effectiveness of the maximal load delivery formulations, a severe contingency analysis use case is considered and a detailed computational study on a variety of realistic network datasets is conducted.  In the interest of making the evaluation exacting, we consider 1000 randomly selected $N$-$k$ scenarios where 30\% of the network's branches are removed.  This type of scenario was observed to be particularly challenging in \cite{COFFRIN2015144}.  The evaluation is presented in four parts: (i) the details of the test cases and computational tools are introduced, (ii) a detailed study of the reliability and scalability of each formulation is conducted, (iii) the optimality gap of the relaxations is analyzed, and (iv) proof-of-concept maximal load delivery analysis is conducted to demonstrate the value of the proposed formulations.

\subsection{Test Cases and Computational Setting}

In this study, seven cases are selected from the IEEE PES PGLib AC-OPF v17.08 benchmark library \cite{pglib_opf}.  The details of these cases are presented in Table \ref{tbl:cases}.  The selection of these cases was designed to (i) provide a representative sample of some of the most realistic datasets available, and (ii) span a wide range of problem sizes, from small (e.g., 73 buses) to realistic (e.g., $>$2000 buses), to highlight the scalablity properties of the proposed formulations.

\begin{table}[h]
\caption{Test Case Size and $N$-$k$ Damage Scenario Overview.}
\label{tbl:cases}
\centering
\begin{tabular}{|l||r|r||r|r|r|r||r|r|r|r|r|r|r|r|r|r|r|r|}
\hline
Test Case & $|N|$ & $|E|$ & $k$ & Scenarios \\
\hline
\hline
IEEE RTS 96 \cite{matpower}\omitit{\cite{780914}} & 73 & 120 & 36 & 1000 \\
\hline
PSERC 240 \cite{6039476} & 240 & 448 & 134 & 1000 \\
\hline
PEGASE 1354 \cite{6488772} & 1354 & 1991 & 597 & 1000 \\
\hline
RTE 1888 \cite{rte_cases} & 1888 & 2531 & 759 & 1000 \\
\hline
Polish 2383wp \cite{pglib_opf} & 2383 & 2896 & 869 & 1000 \\
\hline
Polish 3120sp \cite{pglib_opf} & 3120 & 3693 & 1108 & 1000 \\
\hline
RTE 6468 \cite{rte_cases} & 6468 & 9000 & 2700 & 1000 \\
\hline
\end{tabular}
\end{table}

All of the maximal load delivery formulations were implemented in Julia v0.5 using the optimization modeling layer JuMP.jl v0.18 \cite{DunningHuchetteLubin2017} and PowerModels v0.8 \cite{power_models,pm_pscc}.  The NLP and SOCP formulations were solved with Ipopt \cite{Ipopt} using the HSL MA27 linear algebra solver \cite{hsl_lib}.  The MISOCP formulation was solved using Pajarito.jl v0.5 \cite{LubinYamangilBentVielma2016} using Ipopt and Gurobi v7.0 \cite{gurobi}.  The MINLP formulation was solved using the default version of Bonmin \cite{Bonmin}, which uses open-source solvers.  All of the solvers were configured to stop when the optimality gap was less than $10^{-6}$.  The continuous formulations (i.e., AC-MLD-C, SOC-MLD-C) were given a time limit of 150 seconds, whereas the discrete formulations (i.e., AC-MLD, SOC-MLD) were given a time limit of 1500 seconds.  All of the formulations were evaluated on HPE ProLiant XL170r servers with two Intel 2.10 GHz CPUs and 128 GB of memory.

\subsection{Algorithm Reliability and Scalability}

Table \ref{tbl:conv_stdy} presents the algorithm status and runtime results of all four formulations across all seven test cases.  The status categories indicate the following properties: {\em converged} - the solver algorithm completed normally; {\em time limit} - the solver algorithm reached the prescribed time limit; and {\em error} - the solver algorithm had an unexpected failure, usually due to numerical accuracy issues.  Note that the runtimes are averaged over all of the scenarios for each combination of formulation and status value.

The results are summarized as follows:
\begin{enumerate}[label=(\roman*)]
\item For the smallest test case, IEEE RTS 96, all of the formulations worked effectively.
\item As the test cases increased in size to the 2000--3000 bus range, the reliability and runtime of the AC-MLD formulation became a significant issue.
\item When the test case was very large ($>$6000 buses), the AC-MLD-C relaxation's runtime became a significant limiting factor.
\item Both the SOC-MLD and SOC-MLD-C were reliable and scalable on all of the considered networks; however, SOC-MLD-C was at least 10 times faster than SOC-MLD.
\end{enumerate}
Given the algorithmic success of SOC-MLD/SOC-MLD-C over AC-MLD/AC-MLD-C, the next important topic to investigate is the bounding quality of these convex relaxations. 

\begin{table*}[t]
\small
\centering
\caption{Model Convergence and Runtime Analysis}
\label{tbl:conv_stdy}
\centering
\begin{tabular}{|r||r|r|r|r||r|r|r|r|r|r|r|r|r|r|r|r|r|r|}
\hline
 & \multicolumn{4}{c||}{Solver Status Breakdown} & \multicolumn{4}{c|}{Average Runtime (seconds)} \\
Status & AC-MLD & AC-MLD-C & SOC-MLD & SOC-MLD-C & AC-MLD & AC-MLD-C & SOC-MLD & SOC-MLD-C \\
\hline
\hline
\multicolumn{9}{|c|}{IEEE RTS 96 ($n$=1000)} \\
\hline
converged & 98.60\% & 94.40\% & 100.00\% & 100.00\% & 16.18 & 0.14 & 0.50 & 0.07 \\
\hline
time limit & 1.40\% & -- & -- & -- & 1526.11 & -- & -- & -- \\
\hline
error & -- & 5.60\% & -- & -- & -- & 150.00 & -- & -- \\
\hline
\hline
\multicolumn{9}{|c|}{PSERC 240 ($n$=1000)} \\
\hline
converged & 71.30\% & 98.70\% & 100.00\% & 100.00\% & 18.46 & 1.54 & 4.10 & 0.32 \\
\hline
time limit & 3.70\% & 0.70\% & -- & -- & 1614.58 & 69.82 & -- & -- \\
\hline
error & 25.00\% & 0.60\% & -- & -- & 1500.00 & 20.63 & -- & -- \\
\hline
\hline
\multicolumn{9}{|c|}{PEGASE 1354 ($n$=1000)} \\
\hline
converged & 94.10\% & 100.00\% & 100.00\% & 100.00\% & 221.36 & 5.46 & 32.97 & 2.26 \\
\hline
time limit & 1.70\% & -- & -- & -- & 1598.14 & -- & -- & -- \\
\hline
error & 4.20\% & -- & -- & -- & 1500.00 & -- & -- & -- \\
\hline
\hline
\multicolumn{9}{|c|}{RTE 1888 ($n$=1000)} \\
\hline
converged & 84.30\% & 88.40\% & 100.00\% & 100.00\% & 99.38 & 10.53 & 67.21 & 2.53 \\
\hline
time limit & 0.10\% & 11.50\% & -- & -- & 1655.67 & 151.00 & -- & -- \\
\hline
error & 15.60\% & 0.10\% & -- & -- & 1500.00 & 130.41 & -- & -- \\
\hline
\hline
\multicolumn{9}{|c|}{Polish 2383wp ($n$=1000)} \\
\hline
converged & 86.00\% & 100.00\% & 99.80\% & 100.00\% & 121.58 & 7.64 & 38.43 & 3.04 \\
\hline
time limit & 1.30\% & -- & -- & -- & 1639.76 & -- & -- & -- \\
\hline
error & 12.70\% & -- & 0.20\% & -- & 1500.00 & -- & 990.95 & -- \\
\hline
\hline
\multicolumn{9}{|c|}{Polish 3120sp ($n$=1000)} \\
\hline
converged & 15.70\% & 99.90\% & 99.80\% & 100.00\% & 659.79 & 9.49 & 67.67 & 4.03 \\
\hline
time limit & 74.00\% & 0.10\% & -- & -- & 1531.44 & 150.75 & -- & -- \\
\hline
error & 10.30\% & -- & 0.20\% & -- & 1500.00 & -- & 633.83 & -- \\
\hline
\hline
\multicolumn{9}{|c|}{RTE 6468 ($n$=1000)} \\
\hline
converged & 14.50\% & 35.60\% & 99.60\% & 100.00\% & 865.42 & 57.17 & 648.47 & 12.34 \\
\hline
time limit & 3.10\% & 58.70\% & -- & -- & 1608.88 & 152.38 & -- & -- \\
\hline
error & 82.40\% & 5.70\% & 0.40\% & -- & 1500.00 & 92.89 & 1500.00 & -- \\
\hline
\end{tabular}
\end{table*}

\subsection{Relaxation Strength}

\begin{figure}[t!]
    \begin{center}
    \includegraphics[width=9.0cm]{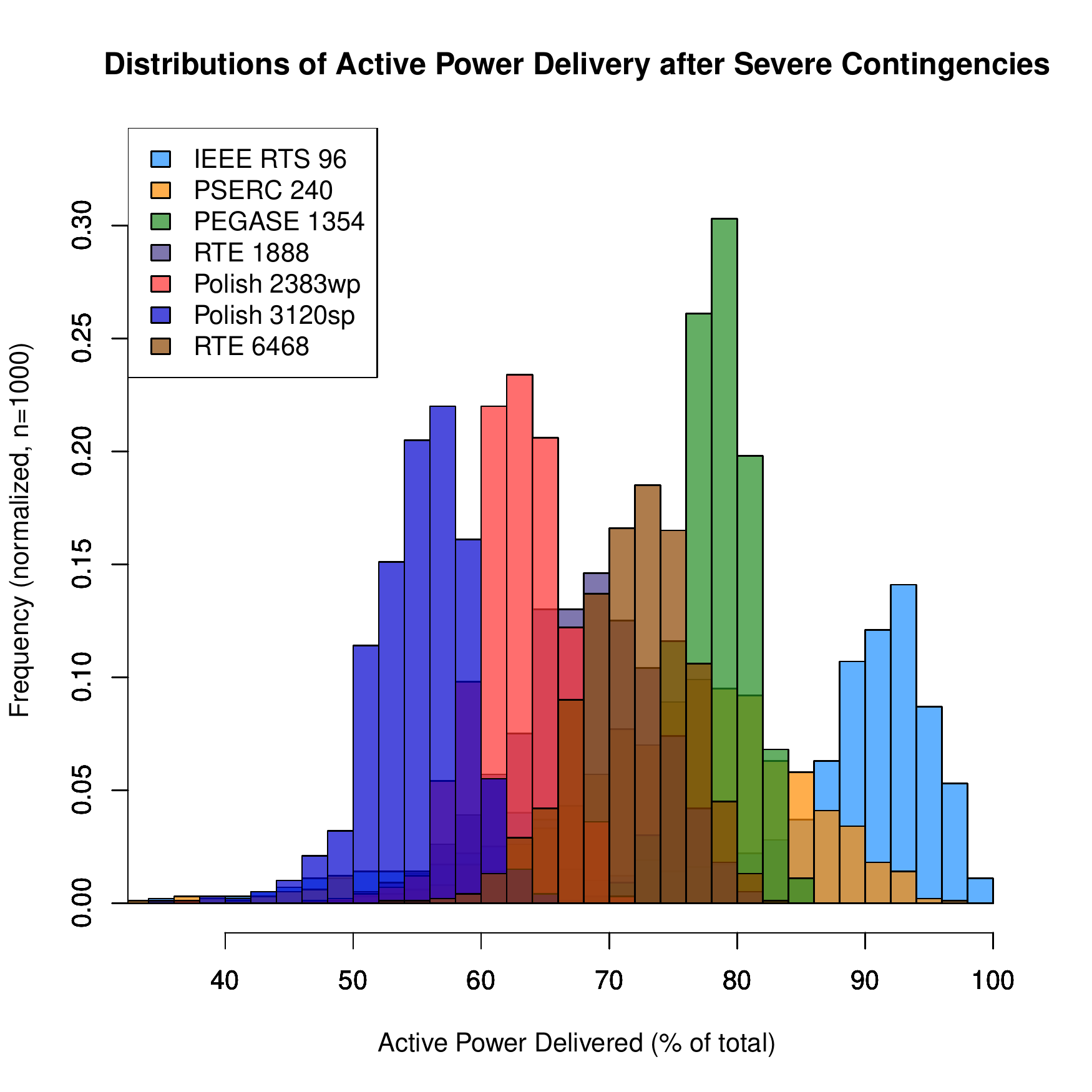}
    \end{center}
    \vspace{-0.6cm}
    \caption{A proof-of-concept maximal load delivery analysis.}
    \vspace{-0.4cm}
    \label{fig:shed_stdy}
\end{figure}

Table \ref{tbl:gap_stdy} presents a comparison of optimality gaps across all four formulations.  Given that the AC-MLD formulation is the only one that provides feasible solutions to this problem, the (percentage) optimality gap is given by
\begin{align}
100 \cdot \left(\frac{\mbox{AC-MLD - Relaxation}}{\mbox{AC-MLD}}\right) \% \nonumber
\end{align}
The values in this table are averages over all scenarios; it is important to note that these averages consider only scenarios where all four algorithms converged. Hence the number of scenarios considered in the average decreases with the network size.

The results are summarized as follows.
All of the optimality gaps were surprisingly small, well below 1\%.
Although there were a few differences in the gaps between SOC-MLD and SOC-MLD-C, the gaps were extremely close.
Overall, these results suggest that the SOC-MLD-C relaxation is sufficient to provide a high-quality upper bound to AC-MLD.  Furthermore, these quality results indicate that the additional runtime needed to solve SOC-MLD is likely not warranted by the small improvement in optimality gap.

\begin{table*}[t!]
\small
\caption{Optimality Gap Study}
\label{tbl:gap_stdy}
\centering
\begin{tabular}{|r||r||r|r|r|r|r|r|r|r|r|r|r|r|r|r|r|r|r|}
\hline
& Obj. Val. & \multicolumn{3}{c|}{Optimality Gap (\%)} \\
     & \multicolumn{1}{c||}{AC-} & \multicolumn{1}{c|}{AC-} & \multicolumn{1}{c|}{SOC-} & \multicolumn{1}{c|}{SOC-}  \\
Case & \multicolumn{1}{c||}{MLS} & \multicolumn{1}{c|}{MLS-C} & \multicolumn{1}{c|}{MLS} & \multicolumn{1}{c|}{MLS-C}  \\
\hline
\hline
IEEE RTS 96 ($n$=914) & 2.463e+04 & 0.0000\% & 0.0044\% & 0.0044\% \\
\hline
PSERC 240 ($n$=707) & 1.945e+06 & -0.0049\% & 0.0010\% & 0.0010\% \\
\hline
PEGASE 1354 ($n$=927) & 2.001e+06 & 0.0022\% & 0.0024\% & 0.0037\% \\
\hline
RTE 1888 ($n$=728) & 1.010e+06 & -0.0267\% & 0.0006\% & 0.0006\% \\
\hline
Polish 2383wp ($n$=859) & 5.759e+05 & 0.0064\% & 0.0062\% & 0.0077\% \\
\hline
Polish 3120sp ($n$=155) & 1.078e+06 & 0.0079\% & 0.0138\% & 0.0186\% \\
\hline
RTE 6468 ($n$=54) & 9.496e+06 & -0.0151\% & 0.0003\% & 0.0003\% \\
\hline
\end{tabular}
\vspace{-0.4cm}
\end{table*}

\subsection{Proof-of-Concept Load Delivery Analysis}

Given that the previous results have demonstrated the efficacy of the SOC-MLD-C relaxation, this section provides a brief and preliminary analysis of the maximal load delivery results provided by that formulation.

Figure \ref{fig:shed_stdy} presents the distribution of active power load delivered in the 1000 $N$-$k$ scenarios for each of the networks considered.  Two observations become clear: (i) in many of these severe contingencies, the amount of load removed can be significant (e.g., $>$20\% of the total active power); and (ii) there do seem to be significant variations in the effects of $N$-30\% contingencies across these networks, in both the mean and variance of total load delivered.  It is also important to recall that these distributions are for the upper bound of active power delivery provided by the SOC-MLD-C relaxation.  The true distributions of AC-MLD may decrease the amount of active power delivered.  However, the results from Table \ref{tbl:gap_stdy} suggest that these changes would be minor.

\section{Conclusion}
\label{sec:conclusion}

This work has proposed a formulation of the AC-MLD problem for use in severe contingency analysis of AC power networks.  The work has shown that even if the complete AC-MLD problem is challenging to compute on large-scale networks (e.g., $>$2000 buses), the convex relaxation SOC-MLD-C is a fast and reliable alternative that provides high-quality bounds for the AC-MLD problem. The value of these high-quality bounds was demonstrated by a proof-of-concept load delivery analysis study.

Despite the notable tightness and fast performance of the SOC-MLD-C relaxation, it is important to emphasize that a tight convex relaxation does not necessarily imply proximity to an AC-feasible solution \cite{coffrin2018convex}.  Hence developing a scalable algorithm for converting the SOC-MLD-C solution into an AC-MLD-feasible solution remains a valuable question for future work.  Other avenues for future work include extension of the AC-MLD formulation and the associated relaxations to multi-timepoint power restoration problems (e.g., \cite{COFFRIN2015144,coffrin2012last,7038388}) that can incorporate cold load pickups, generation ramp rates, and standing phase angle requirements.

\section*{Acknowledgments}
The research presented in this work was supported by the U.S. Department of Homeland Security's National Risk Management Center and by the U.S. Department of Energy through the Los Alamos National Laboratory LDRD program and the Center for Nonlinear Studies.

\bibliographystyle{IEEEtran}

\noindent
LA-UR-18-29573

\end{document}